\documentclass[11 pt,pdftex]{article}
\pdfoutput=1
\usepackage[margin=2.5cm]{geometry}
\usepackage[all]{xy}

\usepackage{authblk, xcolor,graphicx,amsmath,amsthm,amssymb,amscd,enumerate,hyperref} 
\usepackage{tikz}
\usepackage{amsfonts}
\usepackage{pgf} 
\usepackage{mathrsfs}
\usepackage[utf8]{inputenc}

\usepackage{amsmath}
\usepackage{amsthm}
\usepackage{amsfonts}
\usepackage{amssymb}
\usepackage{graphicx}
\usepackage{multicol}
\usepackage[font=footnotesize,labelfont=bf]{caption}
\usepackage{multirow}

\usetikzlibrary{arrows}
\newlength{\matrixheight}

\newtheorem{thm}{Theorem}

\theoremstyle{definition}
\newtheorem{dfn}[thm]{Definition}

\newtheorem{rem}[thm]{Remark}
\newtheorem{exa}[thm]{Example}

\DeclareMathOperator{\gr}{gr}

\DeclareMathOperator{\Gal}{Gal}
\DeclareMathOperator{\disc}{disc}
\DeclareMathOperator{\sign}{sign}

\newcommand{\QQ}{\mathbb{Q}}
\newcommand{\HH}{\mathbb{H}}
\newcommand{\VV}{\mathbb{V}}

\newcommand{\PP}{\mathbb{P}}
\newcommand{\ZZ}{\mathbb{Z}}
\newcommand{\NN}{\mathbb{N}}
\newcommand{\FF}{\mathbb{F}}
\newcommand{\CC}{\mathbb{C}}

\newcommand{\TT}{\mathbb{T}}

\newcommand{\calH}{\mathcal H}

\newcommand{\Mnew}{\mathcal{M}_0}

\usepackage{makecell}
\usepackage{comment}

\usepackage{mathtools}
\usepackage{amsmath}
\usepackage{array}
\usepackage{xparse}
\usepackage{tikz-cd}
\usepackage{tikz}
\usepackage{pgfplots}
\pgfplotsset{compat=1.15}
\usepackage{mathrsfs}
\usetikzlibrary{arrows}
\usepackage[margin=1.5cm]{caption}

\ExplSyntaxOn
\NewDocumentCommand{\sortandprint}{ m }
 {
  \seq_set_split:Nnn \l_tmpa_seq {,} {#1}
  \seq_sort:Nn \l_tmpa_seq
   {
    \fp_compare:nNnTF { ##1 } > { ##2 } { \sort_return_swapped: } { \sort_return_same: }
   }
  \seq_use:Nn \l_tmpa_seq { , }
 }
\ExplSyntaxOff

\usepackage{placeins}

\begin{document}

\title{Hypergeometric local systems over $\QQ$\\ with Hodge vector $(1,1,1,1)$}

\author{Giulia Gugiatti\footnote{Dipartimento di Matematica  ``Tullio Levi-Civita", Universit\`{a}  di Padova, Via Trieste, 63, 35131 Padova.\\  \textcolor{white}{-----} Email address: giulia.gugiatti@math.unipd.it} } 
\author{Fernando Rodriguez Villegas\footnote{Math Section, ICTP, Leonardo Da Vinci Building, Strada Costiera 11, 34151 Trieste.\\ \textcolor{white}{-----} Email address: villegas@ictp.it}}
\affil{}

\date{}

\maketitle

\begin{abstract}
  In this note we consider all irreducible rank-$4$ hypergeometric
  local systems $\HH$ defined over $\QQ$ that support a rational
  one-dimensional variation of Hodge structures of weight $3$ and
  Hodge vector $(1,1,1,1)$. Up to a natural equivalence there are only
  $47$ cases. The first $14$ cases have maximally unipotent monodromy
  at one point and have been extensively studied in the literature. We
  show that all $47$ local systems are associated to families of
  generically smooth threefolds and we analyze the geometry and
  arithmetic at their conifold point.

\end{abstract}

\tableofcontents{}

\section*{Introduction}
\addcontentsline{toc}{section}{Introduction}
In this note we consider all irreducible rank-$4$ hypergeometric local
systems $\HH$ defined over $\QQ$ that support a rational
one-dimensional variation of Hodge structures (VHS) of weight $3$ and
Hodge vector $h=(1,1,1,1)$.  

A hypergeometric local system $\HH$ is
determined by standard hypergeometric parameters $\alpha,\beta$. We
consider them up to isomorphism and pullback by the inversion map
$\iota \colon t \mapsto 1/t$ (which exchanges $\alpha$ and
$\beta$). Here $t$ is our choice of parameter in $\PP^1$ normalized so
that $\HH$ is defined over $\PP^1 \setminus \{0, 1, \infty\}$ and  the local monodromy at $t=1$ is a pseudo-reflection;  following
the Physics literature we call $t=1$ the {\it conifold} of the
family. Up to this equivalence there are only 47 cases. We list them in tables ~\ref{tab:MUM} and ~\ref{tab:NONMUM}.

The first $14$  
cases of our list (table~\ref{tab:MUM})  have maximally unipotent monodromy
(MUM) at one of the singularities (normalized to $t=0$). They have
been extensively studied \cite{BVS, FRVCalabi, DM, 14case} as they correspond
to one-parameter families of Calabi-Yau threefolds and play an
important role in mirror symmetry. We refer to them as MUM cases. The
remaining $33$ cases (table~\ref{tab:NONMUM}) do not have a MUM point.  \smallskip

Our two main goals are:
\begin{enumerate}[1.]
\item
To show  that in fact all of the 47 cases are
associated to families of generically smooth (but non-compact)
threefolds.
\item
To analyze the geometry and arithmetic at the conifold point of each
local system. 
\end{enumerate}
\smallskip

For $\HH$ a hypergeometric local system corresponding to the parameters $\alpha, \beta$, we call the hypergeometric local system $\widetilde{\HH}$
corresponding to the parameters $\widetilde{\alpha},\widetilde{\beta}$
with
$\widetilde{\alpha}_i=\alpha_i+1/2, \widetilde{\beta}_j=\beta_j+1/2$
the \emph{total twist} of 
$\HH$. If $\HH$ is irreducible, the operation of total twist yields an involution. Moreover, if $\HH$ is defined over $\QQ$, then so is $\widetilde{\HH}$, and the rational VHS it supports has the same Hodge vector as that supported by $\HH$.

The 47  cases can be classified as follows:
\begin{itemize}
\item The cases $1, 30, 37$ are equivalent to their total twists.
\item
The cases $2-14$ are the total twists of the cases $15-27$.
\item The cases $28-38$, excluding cases $30$ and $37$, are the total twists
  of the cases $39-47$.
\end{itemize}
Thus our $47$ cases give rise to $25$ pairs
$((\alpha, \beta), (\widetilde{\alpha},\widetilde{\beta}))$.

\paragraph{Goal 1.}
To realize $\HH$ geometrically as a family of threefolds we proceed as 
follows. First, we consider the canonical pair $(Z, \pi)$ of the 
gamma vector $\gamma^{\mathrm{red}}$ corresponding to $\HH$ (see~\S\ref{sec:gamma-vec} and~\S\ref{sec:canonical-toric}). This is typically of high dimension. Then, we apply
to this pair a range of dimension reduction techniques (see~\S~\ref{sec:dim-red}) to
arrive at a family of threefolds. Our techniques in~\S~\ref{sec:gcd-n} extend the dimension reduction method
of~\cite{BCM} (consisting of splitting the 
gamma vector into a union of gamma vectors with relatively prime entries).
We also realize the total twist operation
geometrically in terms of double covers of canonical pairs and fibrations in odd-dimensional quadric bundles~\S~\ref{sec:totaltwist-geometry}. \\
We work out a few cases
explicitly (\S\ref{sec:threefolds}).

\paragraph{Goal 2.}
For each of the $47$ cases $\HH$ degenerates at $t=1$. Local monodromy
at $t=1$ fixes a dimension $3$ subspace with mixed Hodge structure
(MHS) of weights $2$ and $3$ (see \eqref{eq:MHSat1}).  The weight $2$
piece $H_1$ has rank one and Hodge vector $(0,1,0)$, while the weight
$3$ piece $H_2$ has rank $2$ and Hodge vector $(1,0,0,1)$.

Geometrically (\S\ref{sec:t=1-geo}), $H_1$ arises from the two rulings of the even-dimensional quadric
determined by the Hessian $H$ of the polynomial defining the fiber $Z_1=\pi^{-1}(1)$. 
This quadric is the exceptional divisor of the blowup of $Z_1$ at its unique ordinary double
point. On the other hand, $H_2$ is related to
the middle cohomology of the resolution of $Z_1$. 
 
Correspondingly, at the arithmetic level (\S\ref{sec:t=1-ari}), the Euler factor $L_p$ of the
associated $L$-series for good primes $p$~[\S 11-12]\cite{R-V} factors as
\begin{equation}
    \label{eq:euler-factor}
L_p(T)=  \left(1-\sigma(p)p  T\right) \left(1-a_pT+p^3T^2 \right)
\end{equation}
where $\sigma$ is a quadratic character and $a_p$ is
conjecturally\footnote{This is known for the MUM cases~\cite{Schoen,
    FRVCalabi, Y}.}  the $p$-th coefficient of a Hecke eigenform $f$
of weight $4$ and some level. 

We show that $\sigma$ corresponds to the quadratic extension $K/\QQ$
given by the field of definition of the two rulings of the quadric in
the blow-up of $Z_1$. The field $K$ equals $\QQ(\sqrt{\disc(H)})$,
where $\disc(H)$ is the signed discriminant of the Hessian $H$. We
give a formula (equation~\eqref{eq:quad-field}) for $\disc(H)$ in
terms of~$\gamma$.  (These facts about the rank $1$ piece $H_1$ hold
for general symplectic irreducible hypergeometric local systems $\HH$
of any rank defined over $\QQ$.) In the tables we list the
discriminant~$D$ of~$K$ (which equals $\disc(H)$ up to squares).

We identify all modular
forms $f$, excluding those corresponding to cases 38 and 47 (whose	
conductor $5400$ is too large), by computing several $a_p$'s numerically and then searching
in the database LMFDB~\cite{lmfdb}.

Finally, using the determination of $K$, we also identify a rank one
local system $\epsilon_t$, corresponding to a quadratic extension of
the form $\QQ(\sqrt{E t})/\QQ(t)$ for some integer $E$, with the
property that $\widetilde\HH \simeq \epsilon_t\otimes\HH$. We give a
formula for $E$, up to squares, in terms of $\gamma$ (equation~\eqref{eq:quad-field-twist}), and in the tables we also list
$E$.

\paragraph{Related and future work} Detailed statements and complete
proofs of the results in this note, along with explicit equations for
all $47$ families of threefolds, will appear in a subsequent
publication. The MUM cases of our list underlie one-parameter families
of Calabi--Yau threefolds.  The threefolds are mirrors to certain
Calabi--Yau weighted complete intersections $D$, and the families are
(expected to be) Landau--Ginzburg mirrors to certain Fano $4$-fold
weighted complete intersections $X$ with $D$ as anticanonical
section. It would be interesting to understand whether all the 47
cases underlie families of Calabi--Yau threefolds, and/or play a role
in mirror symmetry.

\paragraph{Acknowledgments}
This note expands on the online talk of the second author in occasion
of the series of virtual seminars ``Geometry and Arithmetic of
Calabi-Yau 3-folds", organized by Bhaskaracharya Pratishthana (Pune,
India). He thanks the organizers for the invitation to speak.  We
thank Asem Abdelraouf, Vasily Golyshev, and Vadym Kurylenko for useful
discussions.  We are thankful to Spencer Bloch and Paul Hacking for
helpful correspondence.

\section{Geometric realisation}
\label{sec:geometric-realisation}
In this section we address the problem of the geometric realisation of
the $47$ local systems in tables~\ref{tab:MUM} and~\ref{tab:NONMUM} as families of
threefolds. We start by spelling out our notion of gamma vector and by
reviewing the construction of the canonical pair associated to a gamma
vector with relatively prime entries.
\subsection{Family parameters and gamma vectors}
\label{sec:gamma-vec}
For $\HH$ an irreducible hypergeometric local system of rank $m$ corresponding to parameters $\alpha=(\alpha_0, \dots, \alpha_{m-1}), \beta=(\beta_0, \dots, \beta_{m-1})$, consider the rational function $Q=q_\infty/q_0$, 
where:
\begin{equation}
\label{eq:char-polys}
q_\infty= \prod_{j=0}^{m-1} \left(x-e^{2 \pi \mathrm{i} \alpha_j}\right)  \qquad q_0=\ \prod_{j=0}^{m-1}  \left(x-e^{2 \pi \mathrm{i} \beta_j}\right)  
\end{equation} 
The function $Q$ is defined over the subfield $F \subset \CC$ generated by the coefficients of $q_\infty$ and $q_0$. Following \cite[\S~2]{R-V}, we call $Q$ the \emph{family parameter}, and we 
say that \emph{$\HH$ is defined over $F$}.

If $\alpha_i, \beta_j$ are rational numbers, so that $F$ is a subfield of a cyclotomic field, $\HH$ has a motivic interpretation \cite[$\S5.4$]{Katz}.
In this note we restrict to cases where
$F=\QQ$.  Then $\HH$ supports a rational   VHS  whose Hodge vector can be computed combinatorially \cite{ACVG, Fedorov,FRV-duke}.

An equivalent way to write $Q$ is as follows
\begin{equation}
    \label{eq:Q-gamma}
    Q= \prod_{i=1}^l \left( x^{|\gamma_i|}-1\right)^{-\mathrm{sign}(\gamma_i)}
  \end{equation} for certain non-zero integers $\gamma_i$ with zero
  sum. The order of the $\gamma_i$ is irrelevant for
  \eqref{eq:Q-gamma}, and we may hence normalize
  $\gamma:=(\gamma_1,\ldots,\gamma_l)$ by requiring it to be  weakly-increasing.

  We call a vector $\gamma$ of non-zero integer entries that sum to
  zero and cannot be all grouped in pairs of integers of opposite
  signs a \emph{gamma vector}.  For example, $(-3,2,1)$ is a gamma
  vector, and so is $(-3,-1,2,1,1)$, but $(-2,-1,1,2)$ is not.  (This
  is a minor extension of the notion of gamma vector in \cite{R-V}.)
We say that a gamma vector $\gamma$ is associated to $\HH$ if it
  satisfies \eqref{eq:Q-gamma}.  \smallskip Note that if
  $\gamma=(\gamma_1, \dots, \gamma_l)$ satisfies \eqref{eq:Q-gamma}
  then so does
$\gamma^{(n)}\coloneq (\gamma_1, \dots, \gamma_l,-n,+n)$,
  $n \neq 0 \in \NN$, but there is a unique gamma vector
$\gamma^{\mathrm{red}}$ associated to $\HH$ with no pair of entries
  summing to zero.

If $\gamma$ is a gamma vector for ${\HH}$, the vector $\widetilde{\gamma}$ obtained by replacing each odd entry $\gamma_i$ of $\gamma$ with the pair $2\gamma_i, -\gamma_i$ is a gamma vector for its total twist $\widetilde{\HH}$. Indeed, letting $Q$ and $\widetilde{Q}$ be the correspinding family parameters, one has
\begin{equation}
\label{eq:widetildeQ-Q}
\widetilde{Q}(x)=
Q(-x)=
\prod_{ \gamma_i  \text{ even}} \left(x^{|\gamma_i|} -1\right)^{-\mathrm{sign(\gamma_i)}}  \prod_{\gamma_i  \text{ odd}} \left(\frac{x^{|2\gamma_i|} -1}{x^{|\gamma_i|} -1}\right)^{-\mathrm{sign(\gamma_i)}}  
\end{equation}
We call $\widetilde{\gamma}$ the total twist of $\gamma$. 
In general $\widetilde{\gamma} \neq \widetilde{\gamma}^{\mathrm{red}}$; for example, the total twist of $\gamma=(-6,-1,2,2,3)$ is $\widetilde{\gamma}=(-6,-3,-2,1,2,2,6)$, and $\widetilde{\gamma}^{\mathrm{red}}=(-3,1,2)$.

\smallskip

In tables \ref{tab:MUM} and \ref{tab:NONMUM}, third column, we list the vectors $\gamma^{\mathrm{red}}$ for all the $47$ cases. For reasons of space we write $n^a$ if the entry $n$ appears $a$ times.
Note that for all cases $\gamma^{\mathrm{red}}$ has relatively prime entries 
(this is not true in general).

\subsection{Canonical pairs and toric models} 
\label{sec:canonical-toric}
From now on, $\HH$ is an irreducible hypergeometric local system defined over $\QQ$. 
Let $\gamma=(\gamma_1, \dots, \gamma_l)$ be an  associated gamma vector with \mbox{$\gcd\{\gamma_1, \dots, \gamma_l\}=1$.}
Then a candidate geometric realisation for $\HH$ is the canonical pair $(Z, \pi)$ associated to $\gamma$. 

\begin{dfn}
\label{dfn:canonical-pair}
The \emph{canonical pair}  associated to $\gamma$ is the pair $(Z, \pi)$ given by
\begin{equation}
\label{eq:Z-pi} Z\coloneqq\left(\sum_{i=1}^l z_i=0\right) \subset \TT^{l-1}
\quad \text{and} \quad \pi\coloneqq \frac{z^\gamma }{\Mnew} \colon Z \to \CC^\times 
\end{equation}
 where $\TT^{l-1} \simeq (\CC^\times)^{l-1}$ is a torus with homogeneous coordinates  
$z_j$, $z\coloneqq (z_1, \dots, z_l)$, $z^\gamma\coloneqq \prod_{j=1}^lz_j^{\gamma_j}$, and 
$\Mnew:=\prod_{j=1}^l{\gamma_j}^{\gamma_j}$ (this equals $M$ in~\cite{FRV-duke}
up to a sign).   
\end{dfn}
The dimension of $Z$ is $d\coloneqq l-2$, and we say that $(Z, \pi)$
has dimension $d$.  It is simple to see that
$Z_t\coloneqq \pi^{-1}(t)$ is singular if and only if $t=1$, and $Z_1$
has a unique ordinary double point.
We expect (see \cite[$\S 4$]{R-V})
the following identity of local systems on
$U=\CC^\times \setminus \{1\}$ to hold:
\begin{equation}
    \label{eq:H-local system}
    \HH= \gr^W_k PR^k_c \pi_{U !}\CC
\end{equation}
where $k\coloneq d-1$ is the dimension of $Z_t$, $\pi_U \colon \pi^{-1}(U) \to U$ is the restriction, and $\gr^W_k PR^k_c \pi_{U \; !}\CC$ denotes the local system with stalk at $t$ the top graded piece of the primitive compactly supported cohomology  $PH^k_c(Z_t, \CC)$.  In this note we assume \eqref{eq:H-local system}.
The arithmetic counterpart of \eqref{eq:H-local system} has been
studied in \cite{BCM}, which 
shows that the hypergeometric series  over $\FF_q$
corresponding to $\gamma$ \cite[Def. 1.1, Thm. 1.3]{BCM} equals the
number of $\FF_q$--rational points of a suitable completion of $Z_t$  (modulo certain corrections). 

\medskip

It can be convenient to consider different equations for $(Z, \pi)$, called toric models in \cite[$\S 4$]{R-V}.
 Toric models arise as specialisation of certain multi-parametric Laurent polynomials  first appearing in \cite{GKZ, Bat, Bat-Bor, Stien}. 
We briefly review the construction.

Pick $d$ vectors  in $\ZZ^l$ 
that,
together with $e \coloneqq (1,\dots, 1)$, span $\ker(\gamma\colon \ZZ^l \to \ZZ)$. 
Denote by $M$ the $d\times l$ integer matrix with the $d$ vectors as rows, and call  ${m}_1, \dots, {m}_l$ its columns.  
Let $\TT^d$ be $d$-dimensional torus with coordinates $x_1, \dots, x_d$. 
Let $u$ be a coordinate on $\CC^\times$, 
let 
$u_1, \dots, u_l$ be coordinates on $\CC^l$, and let $\iota \colon \CC^\times \hookrightarrow \CC^l$ be a map such that
$\prod_{j=1}^l u_j^{\gamma_i} \circ \iota = \mathrm{id}$. 
Then, via the change of coordinates 
sending $z_j \mapsto \iota(u)_j x^{m_j}$, 
the pair $(Z, \pi)$ is the same as the pair formed by 
\begin{equation}
  \label{eq:Z-toric} \left(\sum_{j=1}^l \iota(u)_j x^{m_j} =0 \right) \subset \TT^d \times \CC^\times \quad \text{and} \qquad  \frac{\mathrm{pr_2}}{\Mnew}=\frac u{\Mnew} 
  \end{equation} 
where $x$ is the vector with entries $x_i$ and $\mathrm{pr}_2$ is the projection onto the second factor. 
We say that \eqref{eq:Z-toric} is a \emph{toric model} for $(Z, \pi)$.
Concretely, one may choose a vector $\kappa \in \ZZ^l$ such that $\kappa \cdot \gamma=1$ and set $u_j=u^{\kappa_j}$. 
We  call $M$ and $\kappa$ an \emph{$M$-matrix} and a \emph{$\kappa$-vector} for $\gamma$.

\begin{exa}
\label{exa:28}
Consider case $28$ in table \ref{tab:NONMUM}. 
The vector $\gamma^{\mathrm{red}}$
has length $l=10$, thus the  pair $(Z, \pi)$ is a one-parameter family of $7$-folds. 
We construct a toric model for $(Z, \pi)$. 

We can choose $M$ as:
\begin{equation*}
    \label{eq:construct-M}
M= \begin{bmatrix}
0 & -1 & \textcolor{white}{-}1 & 0 & 0 & 0 & 0 & 0 & 0 & \textcolor{white}{-}0 \\
0 & \textcolor{white}{-}0 & -1 & 1 & 0 & 0 & 0 & 0 & 0 & \textcolor{white}{-}0    \\
0 & \textcolor{white}{-}0 & -1  & 0 & 1 & 0 & 0 & 0 & 0 & \textcolor{white}{-}0 \\
0 & \textcolor{white}{-}0 & -1  & 0 & 0 & 1 & 0 & 0 & 0 & \textcolor{white}{-}0\\
0 & \textcolor{white}{-}1 & \textcolor{white}{-}0  & 1 & 0 & 0 & 1 & 0 & 0 & \textcolor{white}{-}0 \\
0 & \textcolor{white}{-}1 & \textcolor{white}{-}0  & 1 & 0 & 0 & 0 & 1 & 0 & \textcolor{white}{-}0   \\
0 & \textcolor{white}{-}1 & \textcolor{white}{-}0  & 1 & 0 & 0 & 0 & 0 & 1 & \textcolor{white}{-}0 \\
0 & \textcolor{white}{-}0 & -1  & 0 & 0 & 0 & 1 & 0 & 0 & -1\\
\end{bmatrix}
\end{equation*}
\end{exa}
and  $\kappa=(0,0,0,0,0,1,0,0,0,0)$. We obtain the hypersurface:
\[
\left(1+ \frac{x_5x_6x_7}{x_1}+\frac{x_1}{x_2x_3x_4x_8}+ x_2x_5x_6x_7+x_3+ux_4+x_5x_8+x_6+x_7+\frac{1}{x_8}=0 \right) \subset \TT^8 \times \CC^\times
\] 
with the projection onto $u/\Mnew$, where $\Mnew=-27/4$.

\subsection{Dimension reduction}
\label{sec:dim-red} 
In general, the VHS supported by $\HH$ has Hodge weight smaller than the dimension of the fibers of any associated canonical pair $(Z, \pi)$, thus it equals $\gr^W_k PR^k_c \pi_{U !}\QQ$ only up to Tate twist.
Then, it is  natural to ask whether $\HH$ also arises from the variation of cohomology of 
a morphism with fibers of dimension equal to the Hodge weight.
For example, our $47$ local systems all support VHS of weight $3$ but most of them are associated to canonical pairs with odd-dimensional fibers of higher dimension. 
We will show
that we may however find corresponding families of affine threefolds in all cases.
\smallskip

In this section we discuss a range of cases where out of a gamma vector of length $l$ one can build a pair $(Y,w)$ of dimension $d-2j<d$, where $d=l-2$. 

If $\gamma$ is union of $s$ gamma vectors $\gamma^i$ with \mbox{$\gcd\{\gamma^i_j\}=1$}, following \cite[$\S 6 $]{BCM} one may build  a pair $(Y,w)$ of dimension $d-2(s-1)$ out of the canonical pairs $(Z_i, \pi_i)$ of the $\gamma^i$.  
In~\S~\ref{sec:gcd=1} we explain the relation between $(Y,w)$ and $(Z, \pi)$. 

In~\S~\ref{sec:gcd-n} we consider the more general case of a vector $\gamma$ with $\gcd\{\gamma_j\}=1$ that is union of $s$ gamma vectors and a pair of the form $(-n,n)$. Out of $\gamma$ we construct a pair $(Y,w)$ of dimension $d-2s$.
By special cases of our construction, we are able to:
\begin{enumerate}[(1)]
    \item Describe the various canonical pairs associated to $\HH$ in terms of each other. 
    \item Generalise the construction of canonical pair~\S\ref{sec:canonical-toric} to an arbitrary  gamma vector. 
    \item  Generalise \cite[\S6]{BCM} to gamma vectors that are union of $s$ arbitrary gamma vectors. 
\end{enumerate}
We give details here for (2) and (3). 

In~\S~\ref{sec:totaltwist-geometry}, we realize the total twist $\widetilde{\HH}$ geometrically in terms of a double cover of the canonical pair of a gamma vector $\gamma$ for $\HH$, and we relate the construction to the canonical pair of the twist $\widetilde{\gamma}$.

\subsubsection{Union of gamma vectors with relatively prime entries} 
\label{sec:gcd=1}
Assume that $\gamma=(\gamma_1, \dots, \gamma_l)$ splits into $s \geq 2$ gamma vectors $\gamma^{{i}}$, ${i}=0, \dots, s-1$, with $\gcd\{\gamma^i_j\}=1$.  
We write
$
\gamma=\gamma^0 \cup \dots \cup \gamma^{s-1}
$ and, for  
convenience, we order the entries of $\gamma$ according to the splitting. We denote by $l_i$ the length of $\gamma^{i}$ and let $d_i\coloneq l_i-2$.  

\paragraph{The pair $(Y,w)$.}
\label{sec:(Y,w)}
The paper \cite[$\S 6 $]{BCM} associates to $\gamma$ a pair $(Y,w)$ of dimension $d-2(s-1)$ as follows. 
For all $i$ let  $(Z_i, \pi_i)$ be the canonical pair associated to $\gamma^{i}$. 
Then $(Y,w)$ is defined by: 
\begin{equation}
\label{eq:(Y,w)}  
Y\coloneqq Z_1 \times \dots \times Z_s \subset \TT^{l_0-1} \times  \dots \times \TT^{l_{s-1}-1} \quad \text{and} \quad w \coloneqq  \pi_1 \cdot \dots \cdot \pi_s \colon  Y \to \CC^\times
\end{equation}
The hypergeometric sum over $\FF_q$ corresponding to $\gamma$ is related to the number of $\FF_q$-rational points of the varieties $Y_t$ (\cite[Thm. 6.1]{BCM}).

\begin{rem}
\label{rem:convolution}
There is a notion of convolution $\star$ of local systems \cite[Ch.~5]{Katz}. 
It is a standard fact that, if a local system $\VV$ is the convolution $\VV=\VV_0 \star \dots \star \VV_{s-1}$  of $s$ local systems $\VV_i$ and $(Z_i, \pi_i)$ is a geometric realisation of $\VV_i$, then $(Z_0 \times \dots  \times Z_s, \pi_i \cdot \dots \cdot \pi_s)$  is a geometric realisation of $\VV$. 

In our case, the fact that $\gamma=\gamma^0 \cup \dots \cup \gamma^{s-1}
$  corresponds to $\HH$ being the convolution $\HH=\HH_0 \star \dots \star \HH_{s-1}$, where $\HH_i$ is the local system corresponding to $\gamma^i$. 
This explains the construction of $(Y,w)$ from a geometric perspective. More details on the VHS associated to $(Y,w)$ are given below.
\end{rem}

\paragraph{Connection with $(Z, \pi)$.}
We pass to toric models for  $(Z_i, \pi_i)$ to make the relation between $(Y,w)$ and  $(Z, \pi)$ manifest. Note that the difference between the dimensions of the two pairs 
is $d-(l-2s)=2(s-1)$. 
\smallskip

It is simple to see that $(Y,w)$ can be rewritten as:
\begin{equation}
\label{eq:Y-w-final}
   Y=(F_0=F_1=\dots =F_{s-1}=0) \subset 
\TT^{l-s-1} \times \CC^\times  \qquad w=\frac{u}{\Mnew}
\end{equation}
where\footnote{We are replacing the $u$-coordinate $u_0$ of a toric model for $(Z_0, \pi_0)$ with the ratio $\frac{u}{\prod_{i=1}^{s-1} u_i}$. Of course this may be done for any of the $u_i$.} \begin{equation*}
    \label{eq:Fi}
    F_0=\sum_{j=1}^{l_0} \left({\frac{u}{\prod_{i=1}^{s-1} u_i}}\right)^{\kappa^0_j} x_{\bf{0}}^{m^0_j}, \qquad F_i=\sum_{j=1}^{l_{i}} u_{i}^{\kappa^{i}_j} x_{\bf{i}}^{m^{i}_j} \quad  i\geq 1
\end{equation*}
and ${x_{\bf{i}}}_j$ ($i=0, \dots, s-1, j=1, \dots, d_i$), $u_1, \dots, u_{s-1}$ are coordinates on $\TT^{d_0+\dots+d_{s-1}} \times (\CC^\times)^{s-1} \simeq 
\TT^{l-s-1}$, while $u$ is a coordinate on $\CC^\times$.

On the other hand, a toric model for $(Z, \pi)$ is given by:
\begin{equation}
\label{eq:Z-wi}
\left( F_0+w_1 F_1+\dots +w_{s-1}F_{s-1}=0 \right) \subset \TT^{d} \times \CC^\times
\qquad \text{with} \quad \frac{\mathrm{pr}_2}{\Mnew}\end{equation} 
where ${x_{\bf{i}}}_j$, 
$u_1, \dots, u_{s-1}$, $w_1, \dots, w_{s-1}$ are coordinates on $\TT^{l-s-1} \times (\CC^\times)^{s-1}  \simeq \TT^{d}$.

\begin{rem}
If a polytope $P$ has width one and projects onto the $(s-1)$-dimensional simplex, any Laurent polynomial $F$  with support $P$ can be written as $F=F_0+w_1F_1+\dots +w_{s-1}F_s$.  This is known as \emph{Cayley trick}. It is simple to see that $\gamma=\gamma^0 \cup \dots \cup \gamma^{s-1}$ if and only if the Newton polytope $P$ of the polynomial $F(\cdot, \Mnew t)$ defining $Z_t$ (see~\eqref{eq:Z-toric}) satisfies these properties.
In \eqref{eq:Z-wi}, in addition, we have chosen the vectors $m_j$ so that each $F_i$ defines a toric model for 
$(Z_i, \pi_i)$. 
\end{rem}

One finds an isomorphism of VHS on  $U=\CC^\times \setminus \{1\}$:
   \begin{equation} 
   \label{eq:VHS-Y-Z}\gr_k^W R^k_c \pi_{U!}\QQ \; (s-1)= \gr_{k-2(s-1)}^W R^{k-2(s-1)}_c w_{U!}\QQ
   \end{equation}
    where $w_U \colon w^{-1}(U) \to U$ is the restriction. (The local system $\HH$ is   \eqref{eq:VHS-Y-Z} tensored with $\CC$.)
    
Indeed, one can consider the partial compactification $Z_t \subset \widehat{Z}_t 
\subset \TT^{l-s-1} \times \PP^{s-1}$, where $\PP^{s-1} \supset \TT^{s-1}$ has local coordinates $w_1, \dots, w_{s-1}$, 
and the first projection $\phi_t \colon Z_t \to \TT=\TT^{l-s-1}$. 
The fibre of $\phi_t$ over a point in $Y_t \subset \TT$ is the projective space $\PP^{s-1}$, while the fiber over a point in $\TT \setminus Y_t$ is a hyperplane in $\PP^{s-1}$. (Note that if $s=2$ $\phi_t \colon Z_t \to \TT$ is the blow-up along $Y_t$.) 
It follows that:
\begin{equation*}
    \label{eq:grk}
    \gr_k^W H^k_c (\widehat{Z}_t, \QQ)= \gr_k^W \left( H^{k-2(s-1)}_c (Y_t, \QQ) (-(s-1))\right)
\end{equation*} 
On the other hand, one can show that the cohomology of $B_t=\widehat{Z}_t\setminus Z_t$ is controlled by those of the spaces $Z_i$ and does not contribute to $\gr_k^W H^k_c(\widehat{Z}_t, \QQ)$.

\subsubsection{Generalisations of $\S\ref{sec:gcd=1}$}
\label{sec:gcd-n}
We generalise~\S~\ref{sec:gcd=1} to a gamma vector $\gamma=(\gamma_1, \dots, \gamma_l)$ with $\gcd\{\gamma_j\}=1$ and of the form $\gamma=\gamma^0 \cup  \dots \cup \gamma^{s-1} \cup (-n,n)$, where $s \geq 1$, $\gcd\{\gamma^i_j\}= D_i \geq 1$, and $n \neq 0 \in \NN$.
We assume for simplicity that $D_0, \dots, D_{s-1}, n$ are pairwise coprime. 

\smallskip

For all $i$, $\frac{\gamma^i}{D_i}\coloneqq (\frac{\gamma^i_1}{D_i}, \dots, \frac{\gamma^i_{l_i}}{D_i})$ is a gamma vector with relatively prime entries.
Let $M_i$ and $\kappa_i$ a $M$-matrix and a $\kappa$-vector for $\frac{\gamma^i}{D_i}$,
pick $a_0, \dots, a_{s-1}, b$ such that $\sum_{i=0}^{s-1} a_i D_i+bn=1$, and set:  
\begin{equation*}
F_0=\sum_{j=1}^{l_0}   \left( \frac{u^{a_0}u_s^n}{\prod_{m=1}^{s-1}u_m^{D_m}} \right)^{k^0_j} x_{\bf{0}}^{m^0_j} \quad F_i=\sum_{j=1}^{l_i}  (u^{a_i}u_i^{D_0})^{k^i_j} x_{\bf{i}}^{m^i_j}\ \quad F_s=u_s^{D_0}+u^b
\end{equation*}
where $i=1, \dots, s-1$. Then, a toric model for the canonical pair $(Z, \pi)$  of $\gamma$ is:
\begin{equation}
    \label{eq:(Z,pi)-general}
\left(F_0+ \sum_{i=1}^{s-1} w_i F_i + w_sF_s=0\right) \subset \TT^{d} \times \CC^\times \quad \text{with} \quad \frac{\mathrm{pr}_2}{\Mnew}
\end{equation}
We define a  pair $(Y,w)$ by:
\begin{equation}
    \label{eq:Y-mostgeneral}
    Y=\left(F_0=F_1=\dots=F_s=0\right) \subset \TT^{d-s} \times \CC^\times  \quad \text{and} \quad w=\frac{u}{\Mnew} \quad 
\end{equation}
Letting $\widehat{Z}_t \subset \TT^{d-s} \times \PP^s$ be the partial compactification of $Z_t$ in the coordinates $w_i$,  one can now show that:
\begin{equation}
\label{eq:Hgeneral}
\gr_k^W H^k_c (\widehat{Z}_t, \QQ)= \gr_k^W \left( H^{k-2s}_c (Y_t, \QQ) (-s)\right) \quad \text{and} \quad \gr_k^W H^k_c ({Z}_t, \QQ) \hookrightarrow \gr_k^W H^k_c (\widehat{Z}_t, \QQ)
\end{equation}

\paragraph{Application (2): non-relatively prime entries.} 
Let $\gamma=(\gamma_1, \dots, \gamma_l)$ be a gamma vector with $\gcd\{\gamma_j\}=D>1$. 

As a special case of~\S\ref{sec:gcd-n} ($s=1, D_0=D,n=1, a_0=0, b=1$), we find that the canonical pair $(Z^{(1)}, \pi^{(1)})$ associated to the vector $(\gamma_1, \dots, \gamma_l, -1,1)$ is the 
hypersurface 
\begin{equation*}
    \label{eq:Z-11}
    \left( \sum_{j=1}^l u_1^{k_j} x^{m_j} +w_1(u_1^D+u)=0\right) \subset \TT^{d+2} \times \CC^\times  \quad \text{with} \ -\frac{\mathrm{pr}_2}{\Mnew}
\end{equation*}
Then, to  $\gamma$ we can associate the $d$-dimensional pair $(Y,w)$: 
\begin{equation}
    \label{eq:Ycgd}
Y=\left( \sum_{j=1}^l u_1^{k_j} x^{m_j} =0\right) \subset \TT^d \times \CC^\times \qquad  \quad \text{and}  \qquad  \quad  w=\frac{u_1^D}{\Mnew}
\end{equation} 
This generalises the construction of canonical pair~\S~\ref{sec:canonical-toric}.

The pair fits into the diagram:
\begin{equation}
 \label{eq:diagram-Yc=gcd}
  \begin{tikzcd}
Z =Y \dar{\pi}    \drar{w}\\
\CC^\times \rar{q_D}& \CC^\times
\end{tikzcd}
\end{equation}
where $(Z, \pi)$ is the canonical pair of $\frac{\gamma}{D}$, and $q_D$ is the $D$:$1$ cover of $\CC^\times$
induced by the map $-u \mapsto u_1^d$. In particular,
\begin{equation}
\label{eq:localsystemYD}
 R^k_c w_{U !}\QQ=   {q_{D}}_\star R^k_c \pi_{U !}\QQ
\end{equation}

Let $\HH_D$ be the local system corresponding to $\gamma$, and $\HH$ that corresponding to $\frac{\gamma}{D}$.
By \eqref{eq:Hgeneral} and studying $\widehat{Z}^{(1)}_t \setminus {Z}^{(1)}_t$, it is not difficult to see that: 
\begin{equation}
\label{eq:case2-seq2}
0 \to \gr^W_{k+2} R^{k+2}_c \pi^{(1)}_{U !}\CC \to  R^k_c w_{U !}\QQ \to  \gr_k \underline{G}^D \to 0
\end{equation}
where $\underline{G}$ is the constant sheaf with stalk $H_c^{d+1}(\TT^d, \QQ)(1) \simeq \QQ^d$.
It follows that
\begin{equation}
\HH_D={q_{D}}_\star \HH
\end{equation}

\paragraph{Application (3): 
union of arbitrary gamma vectors.}
Let $\gamma=\gamma^0 \cup \dots \cup \gamma^{s-1}$ be a gamma vector  with $\gcd\{\gamma_j\}=1$.

As a special case of ~\S\ref{sec:gcd-n}  (remove $(-n,n)$ and the corresponding contributions, i.e. set $u_s=1$, and ignore $w_s$, $F_s$), we find that $Z$ is isomorphic to
the hypersurface
\begin{equation*}
    \label{eq:Z-union-gcd}
\left(  F_0 +\sum_{i=1}^{s-1} w_i F_i=0 \right) \subset \TT^{d} \times \CC^\times
\end{equation*}
and the variety $Y\subset \TT^{l-s-1} \times \CC^\times$ is defined by the equations $F_0=F_1=\dots =F_{s-1}=0$.

It is easy to see that:
\begin{equation}
\label{eq:Yprod-mostgeneral}
(Y,w)=(Y_0 \times \dots \times Y_{s-1}, \prod_{i=0}^{s-1} w_i)
\end{equation}
where $(Y_i, w_i)$ is the pair associated to $\gamma^i$ by \eqref{eq:Ycgd}. This generalises \eqref{eq:(Y,w)}.

Let $(Z_i, \pi_i)$ be the canonical pair corresponding to $\gamma^i/\gcd\{\gamma^i_j\}$.
By \eqref{eq:Hgeneral} and studying $\widehat{Z}_t \setminus Z_t$, one finds
\begin{equation}
    \label{eq:gammas-sequence}
 \HH=\gr^W_{k}R^{k}_c \pi_{U !} \QQ (s-1)  \hookrightarrow \gr^W_{k-2(s-1)}R^{k-2(s-1)}_c w_{U !} \QQ  
\end{equation}
Unlike in~\S\ref{sec:gcd=1}, the cokernel of the morphism \eqref{eq:gammas-sequence} is not always trivial, and is determined by the middle cohomology of the covers of $Z_i$ defined by $F_i=0$.

\subsubsection{The geometry of total twists}
\label{sec:totaltwist-geometry}
Let $\HH$ and $\widetilde{\HH}$ be a hypergeometric local system and its total twist.
Note that 
 $\widetilde{\HH}$ fits into the split exact sequence of local systems on $U=\CC^\times \setminus \{1\}$:
\begin{equation}
\label{eq:totaltwist-sequence}
0 \to \HH \to  p_{U \star} p_U^{\star} \HH \to \widetilde{\HH}  \to 0
\end{equation}  
where  $p \colon \CC^\times \to \CC^\times$ is the double cover $t \mapsto s^2$, and $p_U \colon p^{-1}(U) \to U$ is the restriction. 
Indeed by the projection formula  \[p_{U \star} p_U^{\star} \HH \simeq \HH \otimes p_{U \star} \CC\] and letting
 $\VV$ be the quotient of $p_{U \star} \CC$ by $\CC$   \eqref{eq:totaltwist-sequence}
is the tensor by $\HH$ of 
\[
0 \to \CC \to p_{U \star} \CC \to \VV \to 0
\]  

\paragraph{Double covers.} 
Assume that $\HH$ is irreducible and defined over $\QQ$, let $\gamma$ be an associated gamma vector and
$(Z, \pi)$ be the canonical pair of $\gamma$ (in fact, any geometric realisation of $\HH$ satisfying \eqref{eq:H-local system}). 

We define a pair $(Z^\prime, \pi^\prime)$ via the commutative diagram 
\begin{equation} \label{eq:Zprime-diagram}
 \begin{tikzcd}
 Z^\prime \rar{P}   \drar{\pi^\prime} \dar{\Pi}  &  Z  \dar{\pi}\\
 \CC^\times \rar{p}& \CC^\times
\end{tikzcd}
\end{equation}
where $p$ is the double cover above. Equivalently, $Z^\prime$ is the double cover of $Z$ determined by $\Mnew \cdot s^2=u$ and \mbox{$\pi^\prime \colon Z^\prime \to \CC^\times$} is the projection onto $t$. 

We find:
\begin{equation} \label{eq:Zprime-totaltwist}
 \begin{tikzcd}
 0 \rar  & \gr_k R^k \pi_{U !} \CC_Z \rar  & \gr_k R^k \pi^\prime _{U !} \CC_Z \rar  &  \gr_k R^k \pi_{ U!} \CC_Z  \otimes \VV \rar  &  0\\
 0 \rar  & \HH \rar \uar  &    p_{U \star} p_U^{\star} \HH  \rar \uar &  \widetilde{\HH} \rar \uar &  0\\
\end{tikzcd}
\end{equation}
\vspace{-1cm}

\noindent where $\pi^\prime _{U} \colon (\pi^\prime)^{-1}(U) \to U$ is the restriction and the vertical arrows are inclusions for $k \geq 0$,  isomorphisms if $k >0$. 
Indeed, by base change, the definition of $\pi^\prime$, and the projection formula, we have that  
\[
 R^k \pi^\prime_{!} \CC_{Z^\prime} =  R^k \pi_! \CC_Z \otimes p_\star \CC
\] which, combined with \eqref{eq:totaltwist-sequence}, implies \eqref{eq:Zprime-totaltwist}.

\paragraph{Relation to the canonical pairs of $\widetilde{\HH}$.}
Let $(\widetilde{Z}, \widetilde{\pi})$ be the canonical pair of the total twist $\widetilde{\gamma}$ of $\gamma$. The dimension of $(\widetilde{Z}, \widetilde{\pi})$ is $\widetilde{d}=l_0+d$, where $l_0=2m$ is the number of odd entries of $\gamma$.
Moreover,
$
\widetilde{\mathcal M}_0 =4^{L_o} \Mnew $, where $L_o$ is the sum of the odd entries of $\gamma$.

Let $x_1, \dots, x_d, y_1, \dots,  y_{l_o}$ be coordinates on $\TT^{d} \times \TT^{l_o} \simeq \TT^{\widetilde{d}}$. Let $\widetilde{u}$ be a coordinate on $\CC^\times$ and set 
$u \coloneqq 4^{L_e} \widetilde{u}$ where $L_e$ is the sum of the even entries of $\gamma$.
Let $M=[m_1, \dots, m_l]$ and $\kappa=(\kappa_1, \dots, \kappa_l)$ be a $M$-matrix and a $\kappa$-vector for $\gamma$. 
Then, by the definition of $\widetilde{\gamma}$, a toric model for $(\widetilde{Z}, \widetilde{\pi})$ is 
\begin{equation}
\label{eq:widetildeZ-toric}
 \left(  \sum_{\gamma_j \text{ odd}}  u^{\kappa_j} (2y_j -y_j^2) x^{m_j} + \sum_{ \gamma_j \text{ even}}  u^{\kappa_j}x^{m_j} =0 \right) \subset  \TT^{\widetilde{d}} \times \CC^\times \quad \text{with}  \quad \frac{\mathrm{pr}_2}{\widetilde{\mathcal M}_0}=\frac{\widetilde{u}}{\widetilde{\mathcal M}_0}
\end{equation} 

For $t \neq 1$ let $\widetilde{Z}_t= \widetilde{\pi}^{-1}(t)$ be the fibre .  
The partial compactification $\widetilde{Z}_t \subset X_t \subset \TT^d \times \PP^{l_o}$, where $\PP^{l_o}$ has local coordinates $y_j$, is smooth.
By \eqref{eq:widetildeZ-toric},
the projection $\psi_t \colon X_t \to \TT^d$ is a quadric bundle with $(l_0-1)$-dimensional fibres. 

\begin{rem}
\label{rem:width2}  
If a polytope $P$ has width two and projects onto the $2$-dilation of the $r$-dimensional simplex, then any Laurent polynomial with support $P$ defines a (non-proper) quadric bundle with $(r-1)$-dimensional fibers. By definition of $\widetilde{\gamma}$, the polytope $\widetilde{P}$ associated to the fibre $\widetilde{Z}_t$ satisfies these properties (with $r=l_0$).
\end{rem}

It is not difficult to see that
the discriminant of $\psi_t$ is the fibre $Z_t=\pi^{-1}(t)$. Moreover,  
the double cover $P_t \colon Z^\prime_t \to Z_t$ parametrizing degenerate quadrics is the map $P(\cdot, t)$ in \eqref{eq:Zprime-diagram}.
Indeed, by \eqref{eq:widetildeZ-toric}, one finds  that $P_t$ is determined by the square root of the product
\[
\delta_t=\prod_{\gamma_j \text{ odd}} u^{k_j} x^{m_j}
\]
hence, since  $\kappa \cdot \gamma=1, \sum\gamma_jm_j=0$ imply $\sum_{\gamma_j \text{ odd}} k_j=1, \sum_{\gamma_j \text{ odd}} m_j=0$ mod $2$, 
by the square root of $u$.

One can show that $\gr^W_k H^c(\widetilde{Z}_t, \CC)= \gr^W_k H^c(X_t, \CC)$.
By this fact and the natural generalisation of \cite[Theorem 3]{ACGGFRV} to quadric bundles with odd-dimensional fibres,
one finds the short exact sequence:
\begin{equation}
  0  \to \gr^W_k R^k \pi_{U !} \CC_Z \to \gr^W_k R^k \pi^\prime_{U !}\CC_{Z^\prime}  \to  \gr^W_k R^k \widetilde{\pi}_{U !} \CC_{\widetilde{Z}}  \to 0
\end{equation}
This is the first row of 
\eqref{eq:Zprime-totaltwist}.

\subsection{One-parameter families of threefolds}
\label{sec:threefolds}
By the dimension reduction techniques outlined in~\S\ref{sec:dim-red}, one can find a corresponding family of threefolds for all  $47$ cases of our list. 
Below we discuss a few examples and strategies. 
\medskip

\noindent\textbf{Case 28.}
We continue with example \ref{exa:28}.  The vector $\gamma^{\mathrm{red}}$ has length $10$ and splits into the union $\gamma^\mathrm{red}=\gamma^0 \cup \gamma^1 \cup \gamma^2$, where
\[
\gamma^0\coloneqq(-4,-1,2,3) \quad  \gamma^1\coloneqq(-1,-1,2) \quad \gamma^2\coloneqq(-1,-1,2)
\]
thus we can consider the family of threefolds $(Y,w)$ in \eqref{eq:(Y,w)}. 
Choosing
\[
M_0 =  \begin{bmatrix} 1& 0 & 2 & 0\\
0 & 2 & 1 & 0
\end{bmatrix}
\quad 
M_1 = \begin{bmatrix} 2 & 0 & 1  \end{bmatrix}
\quad  M_2= \begin{bmatrix} 2 & 0 & 1  \end{bmatrix}
\]
and $\kappa^0=(0,-1,0,0)$, $\kappa^1=\kappa^2=(0,-1,0)$, by \eqref{eq:Y-w-final} $(Y,w)$ can be written as
\[
Y=\left(  x_1+ \frac{1}{u}x_2^2x_5x_6+x_1^2x_2+1=x_3^2+\frac{1}{x_5}+x_3=x_4^2+\frac{1}{x_6}+x_4=0 \right) \subset \TT^6 \times \CC^\times 
\] where the $x_j$ are the coordinates on $\TT^6$ and $u$ that on $\CC^\times$,
and $w$ is the projection onto $-4u/27$.  
Note that the fibers $Y_t=w^{-1}(t)$ may be written as hypersurfaces in $\TT^4$.

One may also consider the splitting $\gamma^\mathrm{red}=(-4,2,2) \cup (-1,-1,-1,3) \cup (-1,-1,2)$ and proceed as in~\S\ref{sec:gcd-n}(3), or find a family of threefolds for case 39 (its total twist), to which the methods in ~\S\ref{sec:gcd=1} and~\S\ref{sec:gcd-n}(3) also apply, and use~\S\ref{sec:totaltwist-geometry}.
\medskip

\noindent \textbf{Case 30.} The vector $\gamma^{\mathrm{red}}$ has length $12$ and splits into the union  $\gamma^\mathrm{red}=\gamma^0 \cup \gamma^1 \cup \gamma^2 \cap \gamma^3$, where
\[
\gamma^0\coloneqq(-6,3,3) \quad  \gamma^1=\gamma^2=\gamma^3\coloneqq(-1,-1,2) 
\] thus we can obtain a family of threefolds $(Y,w)$ by~\S\ref{sec:gcd-n}(3). 
Choosing 
\[
M_0=[1,0,2] \quad M_1=M_2=M_3=[2,0,1]
\]
$\kappa^0=(0,1,0), \kappa^i=(0,-1,0), i=1,2,3$, and $\kappa=(0,0,0,1,0,1,0,0,0,0,0,0)$,
$Y$ can be written as:
\[
\left( x_1+\frac{1}{u_1u_2u_3}+x_1^2=ux_2^2+\frac{1}{ u_1^3}+ux_2=x_3^2+\frac{1}{u_2^3}+x_3=x_4^2+\frac{1}{u_3^3}+x_4=0
 \right)\subset \TT^{7} \times \CC^\times
\] with $x_i,u_j$ the coordinates on $\TT^7$ and $u$ that on $\CC^\times$, and $w$ is the projection onto $u$.

Alternatively, one could write $\gamma^\mathrm{red}=(-6,-1,2,2,3) \cup (-1,-1,-1,3) \cup (-1,-1,2)$, and use that $(-6,-1,2,2,3)$ is the total twist of $(-3,1,2)$ combined with~\S\ref{sec:totaltwist-geometry} and remark \ref{rem:convolution}. 
Note that case 30 coincides with its total twist. 
\medskip

\noindent\textbf{Case 41.}  The vector $\gamma^{\mathrm{red}}$ has length $14$. One can check that $\gamma^{\mathrm{red}}$ does not split into the union of $s=5$ gamma vectors, thus we cannot arrive at a family of threefold by the methods in~\S\ref{sec:gcd-n}(3). However, one may observe that $\gamma^{\mathrm{red}}=\gamma^0 \cup \gamma^1 \cup \gamma^2 \cap \gamma^3$, where
\[
\gamma^0\coloneqq(-3,-3,6) \quad \gamma^1\coloneqq(-3,-2,-2,1,6) \quad \gamma^2\coloneqq(-2, 1, 1)  \quad \gamma^3\coloneqq(-2, 1, 1)
\] 
and $\gamma^0$ is the total twist of $(-2,-1,3)$. Thus, by ~\S\ref{sec:gcd-n}(2),~\S\ref{sec:totaltwist-geometry}, and remark \ref{rem:convolution}, we can find a family of threefolds $(Y,w)$. 

Concretely, $(Y,w)=(Y_0 \times Z^\prime_1 \times Z_2 \times Z_3, w_0\times \pi^\prime_1 \times \pi_2 \times \pi_3)$
where  
$(Y_0,w_0)$ and $(Z^\prime_1, \pi^\prime_1)$ can be written as:
\[
\begin{split}
Y_0&=\left(x^2+\frac{1}{u}+x=0\right) \subset \TT^1 \times \CC^\times \quad w_0=u^3/64\\ 
Z^\prime_1&=\left( 1+\frac{y^3}{s^2}+y=0\right)  \subset \TT^1 \times \CC^\times  \quad \ \pi^\prime_1=-4s^2/27
\end{split}
\]
while $(Z_i, \pi_i)$, $i=2,3$ are both isomorphic to:
\[
Z=\left(z+v+z^2=0 \right) \subset \TT^1 \times \CC^\times  \quad   \pi=4v
\]

Alternatively, one can find a family of threefolds for case 31 (its total twist), to which the methods in~\S\ref{sec:gcd-n}(3) apply, and use~\S\ref{sec:totaltwist-geometry}. 

\section{The point $t=1$}
\label{sec:t=1}

Let $\HH$ be any of the $47$ local systems in tables~\ref{tab:MUM}and~\ref{tab:NONMUM}, and
let $\mathcal{H}(t)$ be the corresponding family of motives. In this
section we analyze the geometry and the arithmetic of the motive
$\mathcal{H}(1)$ at the point $t=1$. Our analysis naturally
generalises to any symplectic irreducible hypergeometric local system
$\HH$ defined over $\QQ$. We thank S. Bloch for his help on the middle
cohomology of the fibre $Z_1$ and on the action of Galois on the two
rulings of an even-dimensional quadric.

\subsection{Geometric perspective}
\label{sec:t=1-geo}
The local monodromy of $\HH$ at $t=1$, being a pseudo-reflection of
determinant $1$,
has Jordan form
$$
h_1:=\left( 
\begin{array}{cc}
1&1\\
0&1
\end{array}\right)
\oplus I_2.
$$
We see that the invariant subspace $V\subseteq \CC^4$ of $h_1$ is of
dimension $3$. The monodromy filtration on $V$ is of the form
\begin{equation*}
\label{eq:weight-filtr}
0\subseteq W_2\subseteq W_3=V
\end{equation*}
with $\dim W_2 =1$ and $\dim W_3/W_2 = 2$ 
(see~\cite[II \S2]{Kulikov}).  Hence we have a short exact sequence of MHS:
\begin{equation}
\label{eq:MHSat1}
0 \to \gr_2^W V=W_2 \to V\to \gr_3^W V=W_3/W_2 \to 0 
\end{equation}
\vspace{0.01cm}

We interpret the above sequence geometrically as follows. The
canonical pair $(Z, \pi)$ corresponding to a vector $\gamma$
associated to $\HH$ has even dimension $d\coloneqq 2n$.  Let $F$ be a
Laurent polynomial as in \eqref{eq:Z-toric}, so that
$Z_1=(F(\cdot, \Mnew)=0) \subset \TT^d$.

The exceptional
divisor of the blow-up $\phi$: $\tilde Z_1\rightarrow Z_1$ at the unique ordinary double point $p$ is the smooth
$(2n-2)$-dimensional quadric $Q$ defined by:
\begin{equation}
    \label{eq:exceptional-quadric} \left(x H x^T=0\right) \subset \PP^{2n-1}
\end{equation} where $H$ is the Hessian of $F(\cdot, \Mnew)$ at $p$, and $x=(x_1, \dots, x_{2n})$ is a vector of coordinates on $\PP^{2n-1}$. 

Given a smooth $(2n-2)$-dimensional quadric $Q$ in $\PP^{2n-1}$, it is a standard fact that the middle cohomology $H^{2n-2}(Q)$
is generated by the two rulings of $Q$, that is, the two families $\mathcal{L}_1$ and $\mathcal{L}_2$ of codimension $n$ linear subspaces of $\PP^{2n-1}$ contained in $Q$. %

Writing ${\delta}$ and $\eta$ for the difference and the sum of the two rulings and letting $i \colon Q \hookrightarrow \widetilde{Z}_1$ be the inclusion, one has that $i_\star \delta=0$, while $i^\star H^{2n-2}(\widetilde{Z}_1)$ is generated by $\eta$. 
Then, the  long exact sequence of the pair:
\[
\cdots \to H^{2n-2}(\widetilde{Z}_1) \to H^{2n-2}(Q) \to H^{2n-2}(\widetilde{Z}_1, Q) \to H^{2n-1}(\widetilde{Z}_1) \to H^{2n-1}(Q) \to \cdots
\] yields, via the isomorphism $H^{2n-1}(Z_1) =H^{2n-1}(\widetilde{Z}_1,Q)$, the exact sequence  of MHS:
\begin{equation} 
\label{eq:MHSat1-Z1} 0 \to \QQ(-(n-1))\to H^{2n-1}(Z_1, \QQ) \to H^{2n-1}(\widetilde{Z}_1, \QQ) \to 0
\end{equation} which matches with \eqref{eq:MHSat1} up to a Tate twist by $n-2$.

\subsection{Arithmetic perspective} 
\label{sec:t=1-ari}
From the arithmetic point of view 
the exact sequence \eqref{eq:MHSat1} corresponds to the following
factorization of the Euler factor of the $L-series$ of $\calH(1)$~[\S 11-12]\cite{R-V} for
good primes $p$
\begin{equation}
L_p(T)=  \left(1-\sigma(p)p T\right) \left(1-a_pT+p^3T^2 \right),
\end{equation}
where $\sigma$ is a quadratic character and $a_p$ is (conjecturally)
the $p$-th coefficient of a modular form $f$ of weight $4$. This is equation \eqref{eq:euler-factor} in the Introduction. Note that
$\sigma(p)p$ has weight two. 
\smallskip

The character $\sigma$ corresponds to the quadratic extension $K/\QQ$
determined by the square root of 
\begin{equation}
\label{eq:quad-field}
\sqrt{-(-1)^{\frac{d}{2}} \prod_{i=1}^l \gamma_i}
\end{equation}
where $\gamma=(\gamma_1, \dots, \gamma_l)$ is a gamma vector associated to $\HH$ and $d=l-2$. (Note that two different gamma vectors associated to $\HH$ determine by \eqref{eq:quad-field} the same quadratic field.)
We interpret this in terms of the difference $\delta$ of the two rulings of $Q$. 

Given any irreducible hypergeometric local system $\HH$ over $\QQ$ and any canonical pair $(Z, \pi)$ associated to it, it is not hard to show that 
the Hessian $H$ at $p$ of $F(\cdot, \Mnew)$ has determinant  
$\det(H)=-\prod_{i=1}^l \gamma_i$
up to a square factor. 

For a  symmetric $m \times m$ matrix~$A$  over a field $F$ of characteristic
different from 2 (defining a quadratic form), we define the \emph{signed discriminant}
of $A$  as
\begin{equation*}\label{eq:signed-disc}
\disc(A)\coloneq(-1)^\frac{m(m-1)}{2} \det(A)
\end{equation*} This is standard in Witt's theory
\cite[II \S2]{Lam} and the choice of sign is also that  in Picard--Lefschetz
theory~\cite[\S3.2.1]{Voisin}.

Let now $d=2n$ be even.  Given in general a smooth
$(2n-2)$-dimensional quadric $(xAx^T=0) \subset \PP^{2n-1}$ over $F$,
the field of definition $K$ of the two rulings $\mathcal{L}_1$ and
$\mathcal{L}_2$ of the quadric is given by
\[
K\coloneqq 
F\left(
\sqrt{\disc(A)}\right)=F\left( \sqrt{(-1)^n \det (A) }\right)
\]
where the equality holds since $n =(2n)(2n-1)/2$ mod $2$.
Hence $K$ is also the field of definition of the difference of the
rulings~$\delta$. The non-trivial element of $\Gal(K/F)$ acts by
$\delta\mapsto -\delta$.

In our case $F=\QQ$, and the quadratic
extension corresponding to $\sigma$ \eqref{eq:quad-field} is  $K=\QQ(\sqrt{\disc(H)})/\QQ$, the field of definition
of the difference of the
rulings~$\delta$ on $Q \subset \widetilde{Z}_1$.
It is not hard to verify that
$$
\sign(\disc(H))=(-1)^{\frac{r(r-1)}{2}}
$$
where $r$ equals $\#\{\gamma_i<0\}-\#\{\gamma_i>0\}$ or, equivalently,
the number of $\alpha_i$ or $\beta_i$ which are zero modulo $\ZZ$. In
tables~\ref{tab:MUM} and~\ref{tab:NONMUM}, column 4, we list the discriminants $D$ of
the~$47$ quadratic fields.
\smallskip

We have identified the modular forms $f$,
excluding those corresponding to cases 38 and 47, by computing several $a_p$'s numerically and
then searching in the database LMFDB. They are listed in the 5th column of the tables (using the LMFDB label).
\medskip

Comparing the finite hypergeometric sums associated to the motives $\mathcal{H}(t)$ and the twisted motives $\widetilde{\mathcal{H}}(t)$  
(see~\S\cite[Eq. (10.1)]{R-V}),
we deduce that for all $t$ in $\PP^1$:
\begin{equation}
  \widetilde{\mathcal{H}}(t) \simeq \epsilon_t \otimes \mathcal{H}(t) 
\end{equation} where $\epsilon_t$ is a character depending only on the hypergeometric parameters. Comparing the quadratic fields $K$ corresponding to $\mathcal{H}(1)$
and $\widetilde{\mathcal{H}}(1)$, we deduce that $\epsilon_t$ is the
character corresponding to the quadratic field determined by the
square root:
\begin{equation}
    \label{eq:quad-field-twist}
\sqrt{
(-1)^{\frac{l_o}{2}} \Mnew \cdot t }
\end{equation}
where, for a gamma vector~$\gamma$ for $\mathcal{H}(t)$, $\Mnew$ is as
in Definition \ref{dfn:canonical-pair} and $l_o$ is the number of odd
entries of $\gamma$. (As above, different choices of $\gamma$
determine by \eqref{eq:quad-field-twist} the same quadratic field.)
Note that $\epsilon_t$ is the same for $\HH$ and $\widetilde{\HH}$.
In the last column of the tables we list for each case the
discriminant $E$ of the quadratic field for $t=1$, equal
to $(-1)^{l_o}\Mnew$ up to squares.

\section*{Tables}
\addcontentsline{toc}{section}{Tables}
In the two tables we list the 14 MUM cases and the 33 remaining cases. We tabulate the parameters $\alpha, \beta$, the vector $\gamma^\mathrm{red}$, the discriminant $D$ of the quadratic field associated to $\sigma$, the modular form $f$ as labeled in LMFDB, and the discriminant $E$ of the quadratic field associated to $\epsilon_1$.

\FloatBarrier
\begin{table}[ht!]
\centering
\renewcommand{\arraystretch}{1.3}
\begin{tabular}{|l||l|l|l|l|l|}
\hline
$n^\circ$ & $\alpha$ & $\gamma^{\mathrm{red}}$ & $D$ & $f$ & $E$\\
\hline
1 & $\left(\frac{1}{2}, \frac{1}{2}, \frac{1}{2}, \frac{1}{2}\right)$ & $(-2^4, 1^8)$ & 1 & 8.4.a.a & $\textcolor{white}{-}{1}$ \\ \hline
2 & $\left(\frac{1}{3}, \frac{1}{2}, \frac{1}{2}, \frac{2}{3}\right)$ & $(-3,-2^2, 1^7)$ & 12 & 36.4.a.a & $-3$\\ \hline
3 & $\left(\frac{1}{4}, \frac{1}{2}, \frac{1}{2}, \frac{3}{4}\right)$ & $(-4,-2,1^6)$  & 8 & 16.4.a.a & $-4$ \\ \hline
4 & $\left(\frac{1}{6}, \frac{1}{2}, \frac{1}{2}, \frac{5}{6}\right)$ & $(-6,-2, 1^5,3)$  & 1 & 72.4.a.b & $-3$ \\ \hline
5 & $\left(\frac{1}{3}, \frac{1}{3}, \frac{2}{3}, \frac{2}{3}\right)$ & $(-3^2, 1^6)$&  1 & 27.4.a.a & $\textcolor{white}{-}1$\\ \hline
6 & $\left(\frac{1}{4}, \frac{1}{3}, \frac{2}{3}, \frac{3}{4}\right)$ & $(-4,-3,1^5,2)$ & 24 & 9.4.a.a & $~12$ \\ \hline
7 & $\left(\frac{1}{6}, \frac{1}{3}, \frac{2}{3}, \frac{5}{6}\right)$ & $(-6, 1^4, 2)$ & 12 & 108.4.a.a & $\textcolor{white}{-}1$\\ \hline
8 & $\left(\frac{1}{4}, \frac{1}{4}, \frac{3}{4}, \frac{3}{4}\right)$ & $(-4^2, 1^4, 2^2)$ & 1 & 32.4.a.a & $\textcolor{white}{-}1$\\ \hline
9 & $\left(\frac{1}{6}, \frac{1}{4}, \frac{3}{4}, \frac{5}{6}\right)$ & $(-6,-4,1^3, 2^2,3)$  & 8 & 144.4.a.f & $~12$ \\ \hline
10 & $\left(\frac{1}{6}, \frac{1}{6}, \frac{5}{6}, \frac{5}{6}\right)$ & $(-6^2,1^2,2^2,3^2)$ & 1 & 216.4.a.c & $\textcolor{white}{-}1$ \\ \hline
11 & $\left(\frac{1}{5}, \frac{2}{5}, \frac{3}{5}, \frac{4}{5}\right)$ &   $(-5,1^5)$ & 5 & 25.4.a.b & $\textcolor{white}{-}5$\\ \hline
12 & $\left(\frac{1}{8}, \frac{3}{8}, \frac{5}{8}, \frac{7}{8}\right)$ & $(-8,1^4,4)$ & 8 & 128.4.a.b & $\textcolor{white}{-}1$\\ \hline
13 & $\left(\frac{1}{10}, \frac{3}{10}, \frac{7}{10}, \frac{9}{10}\right)$ & $(-10,1^3, 2,5)$ & 1 & 200.4.a.f & $\textcolor{white}{-}5$\\ \hline
14 & $\left(\frac{1}{12}, \frac{5}{12}, \frac{7}{12}, \frac{11}{12}\right)$ & $(-12,-2,1^4,4,6)$ & 1 & 864.4.a.a & $\textcolor{white}{-}1$\\ \hline
\end{tabular}
\caption{The 14 MUM cases. We omit $\beta=(0,0,0,0)$.}
\label{tab:MUM}
\end{table}

\begin{table}[ht!]
\centering
\renewcommand{\arraystretch}{1.3}
\begin{tabular}{|l||l|l|l|l|l|}
\hline
$n^\circ$ & $\alpha$ & $\gamma^{\mathrm{red}}$ & $D$ & $f$ & $E$\\
\hline
 $15 \ (\widetilde{2})$ & $(0,0,\frac{1}{6},\frac{5}{6})$, $(\frac{1}{2}, \frac{1}{2}, \frac{1}{2}, \frac{1}{2})$ & $(-6, -1^7, 2^5,3)$ & $-4$ & 12.4.a.a &$-3$ \\
\hline
$16 \ (\widetilde{3})$ & $(0,0,\frac{1}{4}, \frac{3}{4})$, $(\frac{1}{2}, \frac{1}{2}, \frac{1}{2}, \frac{1}{2})$ & $(-4,-1^6, 2^5)$  & $-8$ &  8.4.a.a &$-4$\\
\hline
$17 \ (\widetilde{4})$ & $(0,0,\frac{1}{3}, \frac{2}{3} )$, $(\frac{1}{2}, \frac{1}{2}, \frac{1}{2}, \frac{1}{2})$  & $(-3, -1^5, 2^4)$ & $-3$ & 24.4.a.a &  $-3$\\
\hline
$18 \ (\widetilde{5})$ & $(\frac{1}{6}, \frac{1}{6}, \frac{5}{6}, \frac{5}{6})$, $(\frac{1}{2}, \frac{1}{2}, \frac{1}{2}, \frac{1}{2})$ & $(-6^2,-1^6,2^6,3^2 )$  & $\textcolor{white}{-}1$ & 27.4.a.a & $\textcolor{white}{-}1$ \\
\hline
$19 \ (\widetilde{6})$ & $(\frac{1}{4}, \frac{3}{4},  \frac{1}{6}, \frac{5}{6})$, $(\frac{1}{2}, \frac{1}{2}, \frac{1}{2}, \frac{1}{2})$ & $(-6,-4,-1^5,2^6,3)$ & $\textcolor{white}{-}8$ & 144.4.a.d & $~12$ \\
\hline
$20 \ (\widetilde{7})$ & $(\frac{1}{3}, \frac{2}{3}, \frac{1}{6}, \frac{5}{6} )$, $(\frac{1}{2}, \frac{1}{2}, \frac{1}{2}, \frac{1}{2})$ &$(-6,-1^4,2^5)$ & $~12$ & 108.4.a.a & $\textcolor{white}{-}1$\\
\hline
$21 \ (\widetilde{8})$ & $(\frac{1}{4}, \frac{1}{4}, \frac{3}{4}, \frac{3}{4})$, $(\frac{1}{2}, \frac{1}{2}, \frac{1}{2}, \frac{1}{2})$ &$(-4^2,-1^4, 2^6)$ & $\textcolor{white}{-}1$ & 32.4.a.a & $\textcolor{white}{-}1$\\
\hline
$22 \ (\widetilde{9})$ & $(\frac{1}{3}, \frac{2}{3}, \frac{1}{4}, \frac{3}{4})$, $(\frac{1}{2}, \frac{1}{2}, \frac{1}{2}, \frac{1}{2})$&$(-4,-3,-1^3,2^5)$ & $~24$ &  72.4.a.a & $~12$\\
\hline
$23 \ (\widetilde{10})$ & $(\frac{1}{3}, \frac{1}{3}, \frac{2}{3}, \frac{2}{3})$, $(\frac{1}{2}, \frac{1}{2}, \frac{1}{2}, \frac{1}{2})$ &$(-3^2,-1^2,2^4)$ & $\textcolor{white}{-}1$ & 216.4.a.c & $\textcolor{white}{-}1$\\
\hline
$24 \ (\widetilde{11})$ & $(\frac{1}{10}, \frac{3}{10}, \frac{7}{10}, \frac{9}{10})$, $(\frac{1}{2}, \frac{1}{2}, \frac{1}{2}, \frac{1}{2})$ &$(-10,-1^5,2^5,5)$ & $\textcolor{white}{-}1$ & 25.4.a.a & $\textcolor{white}{-}5$\\
\hline
$25 \ (\widetilde{12})$ & $(\frac{1}{8}, \frac{3}{8}, \frac{5}{8}, \frac{7}{8})$, $(\frac{1}{2}, \frac{1}{2}, \frac{1}{2}, \frac{1}{2})$ &$(-8,-1^4,2^4,4)$ & $\textcolor{white}{-}8$ & 128.4.a.b & $\textcolor{white}{-}1$ \\
\hline
$26 \ (\widetilde{13})$ & $(\frac{1}{5}, \frac{2}{5}, \frac{3}{5}, \frac{4}{5})$, $(\frac{1}{2}, \frac{1}{2}, \frac{1}{2}, \frac{1}{2})$ &$(-5,-1^3,2^4)$ & $\textcolor{white}{-}5$ & 200.4.a.e & $\textcolor{white}{-}5$\\
\hline
$27 \ (\widetilde{14})$ & $(\frac{1}{12}, \frac{5}{12}, \frac{7}{12}, \frac{11}{12})$, $(\frac{1}{2}, \frac{1}{2}, \frac{1}{2}, \frac{1}{2})$ & $(-12,-1^4, 2^3,4,6)$ & $\textcolor{white}{-}1$ & 864.4.a.a & $\textcolor{white}{-}1$\\

\hline
    28 & $(0, 0, \frac{1}{4}, \frac{3}{4})$, $(\frac{1}{3}, \frac{1}{2}, \frac{1}{2}, \frac{2}{3})$ & $(-4, -1^5, 2^3, 3)$ & $-24$ &  18.4.a.a & $~12$\\
    \hline
    29 & $(0, 0, \frac{1}{4}, \frac{3}{4})$, $(\frac{1}{3}, \frac{1}{3}, \frac{2}{3}, \frac{2}{3})$ & $(-4, -1^4, 2, 3^2)$ & $-8$  &54.4.a.c & $-4$\\
    \hline
    30 & $(0, 0, \frac{1}{6}, \frac{5}{6})$, $(\frac{1}{3}, \frac{1}{2}, \frac{1}{2}, \frac{2}{3})$ & $(-6, -1^6, 2^3, 3^2)$ & $-3$ & 6.4.a.a & $\textcolor{white}{-}1$\\
    \hline
    31 & $(0, 0, \frac{1}{6}, \frac{5}{6})$, $(\frac{1}{3}, \frac{1}{3}, \frac{2}{3}, \frac{2}{3})$ & $(-6, -1^5, 2, 3^3)$ & $-4$ &  54.4.a.d &  $-3$ \\
    \hline
    32 & $(0, 0, \frac{1}{6}, \frac{5}{6})$, $(\frac{1}{4}, \frac{1}{3}, \frac{2}{3}, \frac{3}{4})$ & $(-6, -1^4,  3^2, 4)$ &  $-24$ &  12.4.a.a & $-4$\\
    \hline
    33 & $(0, 0, \frac{1}{6}, \frac{5}{6})$, $(\frac{1}{4}, \frac{1}{4}, \frac{3}{4}, \frac{3}{4})$ & $(-6, -1^3, -2, 3, 4^2)$  & $-4$ & 96.4.a. & $-3$\\
    \hline
    34 & $(0, 0, \frac{1}{6}, \frac{5}{6})$, $(\frac{1}{5}, \frac{2}{5}, \frac{3}{5}, \frac{4}{5})$ & $(-6, -1^4, 2, 3, 5)$ & $-20$ &150.4.a.h & $-15$ \\
    \hline
    35 & $(\frac{1}{3}, \frac{1}{3}, \frac{2}{3}, \frac{2}{3})$, $(\frac{1}{4}, \frac{1}{4}, \frac{3}{4}, \frac{3}{4})$ & $(-3^2, -2^2, 1^2, 4, 4)$ &  $\textcolor{white}{-}1$ & 864.4.a.a & $\textcolor{white}{-}1$\\
    \hline
    36 & $(\frac{1}{3}, \frac{1}{3}, \frac{2}{3}, \frac{2}{3})$, $(\frac{1}{6}, \frac{1}{4}, \frac{3}{4}, \frac{5}{6})$ & $(-3^3,  -2^2, 1^3,  4, 6)$  & $\textcolor{white}{-}8$ &432.4.a.k & $~12$\\
    \hline
    37 & $(\frac{1}{3}, \frac{1}{3}, \frac{2}{3}, \frac{2}{3})$, $(\frac{1}{6}, \frac{1}{6}, \frac{5}{6}, \frac{5}{6})$ & $(-3^4,-2^2, 1^4, 6^2)$ & $\textcolor{white}{-}1$ & 72.4.a.b & $\textcolor{white}{-}1$\\
    \hline
    38 & $(\frac{1}{3}, \frac{1}{3}, \frac{2}{3}, \frac{2}{3})$, $(\frac{1}{10}, \frac{3}{10}, \frac{7}{10}, \frac{9}{10})$ & $(-5, -3^2, -2, 1^3,  10)$ &  $\textcolor{white}{-}1$ & 5400.4.?.? & $\textcolor{white}{-}5$\\
    \hline
   
$ 39 \  (\widetilde{28})$&   $(\frac{1}{4}, \frac{1}{2}, \frac{1}{2}, \frac{3}{4})$, $(0, 0, \frac{1}{6}, \frac{5}{6})$ & $(-4, -3, -2^2,1^5,6)$ & $-8$ & 48.4.a.c &$~12$\\
    \hline
$ 40 \ (\widetilde{29})$ &   $(\frac{1}{4}, \frac{1}{2}, \frac{1}{2}, \frac{3}{4})$, $(\frac{1}{6}, \frac{1}{6}, \frac{5}{6}, \frac{5}{6})$ & $(-4, -3^2, -2^3, 1^4, 6^2)$ & $\textcolor{white}{-}8$ & 432.4.a.j & $-4$\\
    \hline
 $41 \ (\widetilde{31})$ &   $(\frac{1}{3}, \frac{1}{2}, \frac{1}{2}, \frac{2}{3})$, $(\frac{1}{6}, \frac{1}{6}, \frac{5}{6}, \frac{5}{6})$ &  $(-2^4, -3^3, 1^5, 6^2)$ & $~12$  & 54.4.a.a &$-3$\\
    \hline
 $42 \ (\widetilde{32})$ &   $(\frac{1}{3}, \frac{1}{2}, \frac{1}{2}, \frac{2}{3})$, $(\frac{1}{6}, \frac{1}{4}, \frac{3}{4}, \frac{5}{6})$ & $(-2^4, -3^2, 1^4, 4,6)$ &  $~24$   & 48.4.a.a &$-4$\\
    \hline

  $43 \ (\widetilde{33})$&   $(\frac{1}{3}, \frac{1}{2}, \frac{1}{2}, \frac{2}{3})$, $(\frac{1}{4}, \frac{1}{4}, \frac{3}{4}, \frac{3}{4})$ &$(-3,-2^4, 1^3, 4^2)$& $~12$ &288.4.a.f&$-3$\\
    \hline
 $44  \ (\widetilde{34})$&   $(\frac{1}{3}, \frac{1}{2}, \frac{1}{2}, \frac{2}{3})$, $(\frac{1}{10}, \frac{3}{10}, \frac{7}{10}, \frac{9}{10})$ & $( -5,-3,-2^3, 1^4, 10)$ & $~12$ & 450.4.a.o &$-15$\\
    \hline
$45  \ (\widetilde{35})$ &    $(\frac{1}{6}, \frac{1}{6}, \frac{5}{6}, \frac{5}{6})$, $(\frac{1}{4}, \frac{1}{4}, \frac{3}{4}, \frac{3}{4})$ &$(-6^2, -1^2, 3^2,  4^2)$ &$\textcolor{white}{-}1$ &  864.4.a.a & $\textcolor{white}{-}1$\\
    \hline
 $ 46 \ (\widetilde{36})$&   $(\frac{1}{6}, \frac{1}{6}, \frac{5}{6}, \frac{5}{6})$, $(\frac{1}{4}, \frac{1}{3}, \frac{2}{3}, \frac{3}{4})$ &$(-6^2, -1^3,  2, 3^3, 4)$ & $~24$ & 216.4.a.a &$~12$\\
    \hline
    $ 47  \ (\widetilde{38})$ &     $(\frac{1}{6}, \frac{1}{6}, \frac{5}{6}, \frac{5}{6})$, $(\frac{1}{5}, \frac{2}{5}, \frac{3}{5}, \frac{4}{5})$ & $(-6^2, -1^3, 2^2, 3^2, 5)$ & $\textcolor{white}{-}5$ &5400.4.?.? &$~5$\\
    \hline
\end{tabular}

\caption{The 33 cases without a MUM point. We write $n~(\widetilde{m})$ to indicate that the case $n$ is the total twist of a previous case $m$.}
\label{tab:NONMUM}
\end{table}
\FloatBarrier

\footnotesize
\bibliographystyle{amsalpha}
\bibliography{bib-ggfrv}

\end{document}